\documentclass{article}

\usepackage{arxiv}
\usepackage{amsmath}
\usepackage{graphicx}

\usepackage{mathptmx}

\usepackage[utf8]{inputenc} 
\usepackage[T1]{fontenc}    
\usepackage{hyperref}       
\usepackage{url}            
\usepackage{booktabs}       
\usepackage{amsfonts}       
\usepackage{nicefrac}       
\usepackage{microtype}      
\usepackage{graphicx}
\usepackage{doi}

\title{\large Stability of a Parametrically Driven, Coupled Oscillator System: An Auxillary Function Method Approach}

\author{{\includegraphics[scale=0.03]{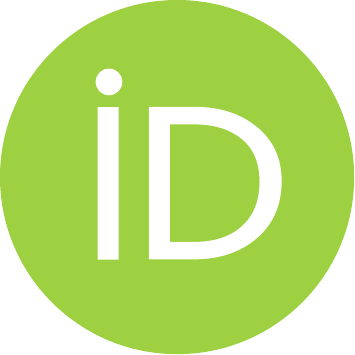}\hspace{1mm}Andrew N. McMillan} \\
	Department of Mathematics\\
	University of Michigan, Ann Arbor\\
	Ann Arbor, MI 48108 \\
	\texttt{andrewnm@umich.edu} \\
	\And
	{\includegraphics[scale=0.03]{orcid.pdf}\hspace{1mm}Yin Lu Young} \\
    Department of Mechanical Engineering and\\ 
    Department of Naval Architecture and Marine Engineering\\ University of Michigan, Ann Arbor\\
    Ann Arbor, MI 48108 \\
	\texttt{ylyoung@umich.edu} \\
	\And
	{\includegraphics[scale=0.03]{orcid.pdf}\hspace{1mm}Mary Robinson} \\
    Department of Physics\\ 
    University of Michigan, Ann Arbor\\
    Ann Arbor, MI 48108 \\
	\texttt{ylyoung@umich.edu} \\
}

\renewcommand{\shorttitle}{Stability of Coupled Oscillator System: An Auxillary Function Method}

\hypersetup{
pdftitle={Stability of Coupled Oscillator System: An Auxillary Function Method},
pdfsubject={q-bio.NC, q-bio.QM},
pdfauthor={Andrew McMillan, Yin Lu Young, Mary Robinson},
pdfkeywords={dynamical systems, Second keyword, More},
}

\begin{document}
\maketitle

\date{\today}
\begin{abstract}
Coupled, nonlinear oscillators are often studied in applied biology, physics, fluids, and many other disciplines. In this paper, we study a parametrically driven, coupled oscillator system where the individual oscillators are subjected to varying frequency and phase with a focus on the influence of the damping and coupling parameters away from parametric resonance frequencies. In particular, we study the key long-term statistics of the oscillator system's trajectories and stability. We present a novel, robust and computationally efficient method come to be known as an auxillary function method for long-time averages, and we pair this method with classical, perturbative-asymptotic analysis to corroborate the results of this auxillary function method. These paired methods are then used to compute the regions of stability for a coupled oscillator system. The objective is to explore the influence of higher order, coupling effects on the stability boundary across a broad range of modulation frequencies, including frequencies away from parametric resonances. We show that both simplified and more general asymptotic methods can be dangerously un-conservative in predicting the true regions of stability due to high order effects caused by coupling parameters. The differences between the true stability boundary and the approximate stability boundary can occur at physically relevant parameter values in regions away from parametric resonance. The differences between the solutions depends on the specific parameters of the system, as explained in the results section. As an alternative to asymptotic methods, we show that the auxillary function method for long-time averages is an efficient and robust means of computing true regions of stability across all possible initial conditions.
\end{abstract}

\maketitle


\section{Introduction}

\par The study of coupled oscillators is ubiquitous in both applied and theoretical fields for their wide range of modeling capabilities \cite{Csaba, Wilson21, Schwemmer2012}, as well as their rich mathematical structure \cite{Fujisaka83, PhysRevE.50.1874}. In applied fields, qualitative analysis of an oscillator system's trajectories or long term statistics in the presence of varying parameters or initial conditions is of primary importance \cite{1075948, PhysRevLett.81.5318}. The stability of coupled oscillators are also of great interest in many fields including engineering, quantum optics, bio-chemistry, etc \cite{Strogatz91, PhysRevE.61.5080, Ermentrout85}. Moreover, in theoretical work, the presence of parametric terms prohibits classical, linear analysis, and novel analysis tools are an interest of development. In this paper, we present a novel, theoretical method to study the long term statistics of a parametrically driven, coupled oscillator system. 
\par To motivate this study, we provide several, key applications of coupled oscillators. In the biological sciences, coupled oscillators have extensively been employed in the study of circadian rhythms where individual circadian oscillations are generated in suprachiasmatic nuclei by an internal regulatory network \cite{Reppert02, Gonze05}. A general interest is to establish both when and how these individual oscillations become synchronized or out of synchrony \cite{Komin10}, and key quantities, such as measures of synchrony, are realized as time averages of the underlying oscillator's correlations \cite{Gonze05, Komin10}. Hence, computing parameter and initial condition dependent time averages is of pivotal importance and frequently challenging. Coupled oscillators have also been explored in the neuroscience and neural network communities, respectively. In neuroscience, for example, the Kuromoto model has received widespread attention. Past and on-going work focuses on the study of long time average frequencies between individual oscillators and investigates how these individual oscillators become globally phase locked \cite{PhysRevLett.109.034103, 2000PhyD..143....1S}. Furthermore, in the study of neural networks, data-driven models attempt to construct potentially nonlinear dynamical models via data sampling \cite{Corbetta20, YEO2019671, Champion22445}. This requires a potentially large data set for the model to be trained on and one's choice of data and assumptions about the form of the model have far reaching consequences on the stability of the underlying dynamics and the accuracy of model predictions. Additionally, coupled oscillators have received widespread attention in the work on fluids. For instance, extensive work has been done on the phenomena of phase lock-ins for vortex induced vibrations \cite{DELANGRE2006783}. The general interest is to establish which range of frequencies induce unstable and growing amplitudes for the underlying dynamics. Finally, there are a variety of applications in physics.  Networks of optical parametric oscillators are being investigated as a way to simulate classical Ising (i.e. spin-$\frac{1}{2}$) systems. Such simulators have exciting potential applications in quantum computation \cite{Marcello19}.  In addition, the behavior of coupled parametric oscillators have been used to study the response of microelectromechanical systems (MEMs) \cite{MEMs}. These examples provide only a small glimpse of the wide range of applications that these coupled oscillator systems have across a variety of scientific disciplines. 
\par Frequently, a primary interest is in studying the effects of the oscillators coupling strength, the type of coupling, the phase differences between the oscillators, and the magnitude, frequency and phase of the driving force on the long term statistics of the system’s solutions. Although most real world applications involve oscillators with differing phases or differing frequencies, most of the analysis that appears in the literature has ignored both of these essential features while also ignoring the effects of coupling strength, especially at frequencies away from parametric resonance \cite{Marcello19}. Ignoring frequencies away from parametric resonance frequencies is natural if one is interested in studying phase-locking or if one has energy harvesting applications in mind. However, real-world systems typically operate over a broad range of parametric, modulation frequencies. Moreover, many engineering and biological applications including coupled oscillator systems are designed to operate away from parametric resonance to avoid excessive vibration, noise, and accelerate fatigue. We demonstrate and explain in Section IV the importance of considering higher order effects across parameter space for the stability boundaries and long-time statistics particularly if the interest is in regions away from parametric resonance frequencies.
\par The reason these essential features have primarily been ignored is almost exclusively due to the problem's difficulty and the limitations of perturbative asymptotic methods \cite{PhysRevLett.69.2472}. From a computational perspective, direct numerical simulation (DNS) is incredibly expensive. Moreover, one is limited to the study of trajectories that are numerically stable, which pivotally depends on one's choice of initial conditions for nonlinear oscillators. Frequently, one is interested in unstable solutions in particular for the ``control of chaos'' \cite{OGY1990}. These concerns naturally pose the problems of which initial conditions to choose and a robust method to determine whether one's findings are generically sensitive to initial data. From a theoretical perspective, analysis is notoriously difficult due to the presence of parametric, nonlinear terms, particularly for coupled systems. Research is frequently concerned with knowing the boundaries of parametric instabilities, especially when the parametric term has a frequency that is different than a multiple of the system's natural frequency. Therefore, approximate solutions are frequently derived via perturbative, asymptotic methods \cite{YANG2015499, CASWELL1979153}. However, these asymptotic methods both fail to capture the full range of the parameters in question and the effects of the oscillator's coupling term, which can make the higher order terms non-negligible. As high order terms are frequently neglected, asymptotics foregoes when or how these higher order terms become relevant \cite{10.2307/4153181}.
\par In this paper, we demonstrate that our proposed method, come to be known as an auxillary function method for long-time averages, is incredibly robust in addressing the difficulties previously explicated. On computational grounds, our method is computationally efficient and foregoes the need to manually check various initial conditions. In fact, our method determines extremal time averages across \textit{all} initial conditions. On theoretical grounds, our method is robust even in the face of highly non-linear dynamics and does not require approximate expansions or restrictions to only small parameter values. Therefore, we can recover the stability boundaries and long term statistics without needing to ignore the higher order, non-linear terms that are disregarded in perturbative asymptotic analysis. Hence, the use of the auxillary function method for stability analysis of nonlinear coupled oscillators is the primary and novel contribution of our work. 
\par Our objective is to study a parametrically driven, coupled oscillator system with the auxillary function method for long-time averages. We demonstrate the robustness of our auxillary function method in computing essential long time averages of the system's trajectories and determine stability boundaries for varying system parameters. We will also compare our results against the results predicted by perturbative asymptotic analysis and investigate the effects of higher order, coupling terms on the system's stability and long-term statistics across a broad range of physically relevant modulation frequencies at and away from parametric resonance. 
\section{Auxillary Function Method for Long-time Averages}
To introduce the auxiliary function method for long-time averages, we focus on determining upper bounds for time averages of functions of the dynamical variable for  autonomous ordinary differential equations (ODEs).
\par Consider $\textbf{x}(t)\in \mathbb{R}^d$ satisfying
\begin{equation}
        \dot{\textbf{x}}=\textbf{f}(\textbf{x})
        \label{eq1}
\end{equation}
for $\textbf{f}:\mathbb{R}^{d}\rightarrow \mathbb{R}^d$ a continuously differentiable vector field. When there is no confusion, we denote $\textbf{x}_i(t)$ and $\textbf{f}_i(\textbf{x})$ as the i-th component of the vectors $\textbf{x}$ and $\textbf{f}(\textbf{x})$, respectively. \par For the purposes of the computation of long term statistics, it becomes natural to ask which trajectories of Equation (\ref{eq1}) realize the optimal time average and their corresponding averaged values.  For a general quantity of interest $\Phi(\textbf{x})$, we define its long-time average along the trajectories $\textbf{x}(t)$ with $\textbf{x}(0) = \textbf{x}_0$ by:
\begin{equation}
    \overline{\Phi}(\textbf{x}_0):=\limsup_{T\rightarrow +\infty}\frac{1}{T}\int_{0}^T\Phi(\textbf{x}(t)) dt.
\end{equation}
As the underlying dynamics are of physical relevance, it is natural to assume that all trajectories $\textbf{x}(t)$ lie in a compact subset B of the phase space $\mathbb{R}^{d}$.
We are interested in the \textit{maximal long-time average} among all trajectories (eventually) remaining in $\text{B}$, i.e.,
\begin{equation}
  \overline{\Phi}^*=\max_{\textbf{x}_0\in \text{B}}\overline{\Phi}(\textbf{x}_0).
  \label{phibarstar}
\end{equation}
Upper bounds on averages can be deduced using the fact that time derivatives of bounded functions average to zero. This elementary observation follows from the fact that for every $\text{V}(\textbf{x}) \in C^1(\text{B})$---the set of continuously differentiable functions on B---we have
\begin{equation}
\begin{split}
    0&=\limsup\limits_{T \rightarrow +\infty}\frac{\text{V}(\textbf{x}(T))-\text{V}(\textbf{x}(0))}{T}\\
    &=\overline{\frac{d}{dt}\text{V}(\textbf{x}(\cdot))}\\
    &=\overline{\textbf{f}(\textbf{x}(\cdot))\cdot\nabla \text{V}(\textbf{x}(\cdot))}.
    \end{split}
\end{equation}
We hereafter refer to any such $\text{V}(\textbf{x}) \in C^1(\text{B})$ as ``auxiliary'' functions.
Note that Equation (4) holds for any auxiliary function so there is an infinite family of functions with the same time average as $\Phi(\textbf{x})$. In particular,
\begin{equation}
    \overline{\Phi}(\textbf{x}_0)=\overline{  \Phi(\textbf{x}_0)+\textbf{f}(\textbf{x}_0) \cdot\nabla \text{V}(\textbf{x}_0)}.
\end{equation}
For any auxiliary function one obtains an upper-bound on $\overline{{\Phi}}^*$ by bounding the right hand-side point-wise on $\text{B}$ and subsequently maximizing the left hand side over initial data $\textbf{x}_0$:
\begin{equation}
\overline{{\Phi}}^* \leq \max_{\textbf{x}\in \text{B}}[\Phi(\textbf{x})+\textbf{f}(\textbf{x})\cdot\nabla \text{V}(\textbf{x})].
     \label{U_upper}
     \end{equation}
The sharpest upper bound on $\overline{{\Phi}}^*$ is then
\begin{equation}
    \overline{{\Phi}}^*\leq \inf_{\text{V}\in C^1(\text{B})}\max_{\textbf{x}\in \text{B}} \ [\Phi(\textbf{x})+\textbf{f}(\textbf{x})\cdot \nabla \text{V}(\textbf{x})],
    \label{ineq}
\end{equation}
and we note that the minimization on the right hand side is convex in the auxillary function $V$. The remarkable fact is that the inequality in Equation (\ref{ineq}) can actually be made an equality. The details of the proof can be found in \cite{Tobasco18}, but the key ideas are that maximizing time averages can be realized as  maximizing phase space averages against invariant measures via Birkhoff's Ergodic Theorem; maximizing phase space averages over invariant probability measures can be realized as a Lagrange multiplier problem, where the auxillary function $V$ can be interpreted as a Lagrange multiplier; finally, the resulting max-min problem can be written as a min-max problem using a variety of abstract min-max theorems over infinite dimensional spaces. 
\par Note that the bound $\overline{\Phi}^*\leq U$ follows from Equation (\ref{U_upper}) if
\begin{equation}
    U-\Phi(\textbf{x})-\textbf{f}(\textbf{x})\cdot \nabla V \geq 0\,\forall \textbf{x}\in \text{B}, 
    \label{nonneg}
\end{equation} 
and if one restricts to polynomial dynamics, then the problem reduces to a convex optimization problem subject to a local, non-negative polynomial constraint.
\par When B is a semi-algebraic set \footnote{A subset B in $\mathbb{R}^d$ is called semi-algebraic if it is defined by a finite collection of polynomial equations of the form $P(\textbf{x})=0$ and $Q(\textbf{x})>0$ for $\textbf{x}\in \mathbb{R}^d$.}, one can augment the constraint in Equation (\ref{nonneg}) with polynomials that define B; this is called the S-procedure. The non-negative polynomial constraint in Equation (\ref{nonneg}) can then be replaced with a sum of squares\footnote{For $\textbf{x}\in \mathbb{R}^d$, a polynomial $p(\textbf{x})$ is said to be a sum of squares if $p(\textbf{x})=\Sigma p_i(\textbf{x})^2$ for some finite collection of polynomials. The set of all sum of squares polynomials over $\textbf{x}\in \mathbb{R}^d$ will be denoted $\mathcal{S}[\textbf{x}]$.} constraint, which has an equivalent semi-definite program formulation. In the resulting semi-definite program, the polynomial degree of the auxillary function $V$ is a degree of freedom, where increasing the degree of $V$ improves the sharpness of the computed bounds. When increasing the degree of $V$ no longer improves the bound, we say the bounds are sharp, at least to numerical precision; see \cite{Tobasco18} for discussions and reviews on sum of squares programs via semi-definite programming. 
\par One should note that the above formulation is completely dependent on the original dynamical system in Equation (\ref{eq1}) being autonomous. For our purposes, the model in question will have non-linear, non-autonomous trigonometric dependence. Therefore, the dynamics generically take on the form: 
\begin{equation}
    \dot{\textbf{x}}=\textbf{f}(\textbf{x},\cos(\omega t),\sin(\omega t)).
    \label{nonauto}
\end{equation}
 The canonical way of making systems of the form seen in Equation (\ref{nonauto}) autonomous is to introduce an additional coordinate $x_{d+1}= t$ and extend the system dimension from $d$ to $d+1$. 
 However, this method unfortunately introduces an unbounded dependent variable while retaining non-polynomial dependence on it. This problem can be circumvented by introducing \textit{two} new dynamical variables satisfying the polynomial sub-system
\begin{equation}
    \begin{split}
        \dot{x}_{d+1}=&-\omega x_{d+2}\\
       \dot{x}_{d+2}=&\omega x_{d+1}\\
        \text{Subject to:} \quad &x^2_{d+1}+x^2_{d+2}=1.
    \end{split}
    \label{4.2}
\end{equation}
One should note that equations of the form in Equation (\ref{4.2}) are frequently referred to as differential-algebraic systems- that is, a differential equation subject to an algebraic constraint. In order to write the system in Equation (\ref{4.2}) as a standard dynamical system, the form seen in Equation (\ref{nonauto}), one needs to omit the origin as a potential solution either via assumption or other techniques.
\par Additionally, one can also formulate equivalent autonomous polynomial dynamics for both quasiperiodic and substantially more complex $\frac{2 \pi}{\omega}$-periodic time dependences in the vector field. Employing a new pair of dynamical variables like those in Equation (\ref{4.2}) for each independent frequency allows for quasiperiodic time dependence, at least for quasiperiodicity involving only a {\it finite} number of independent frequencies. Other $\frac{2 \pi}{\omega}$-periodic time functions can be expressed  as finite linear combinations of $\cos(n \omega t)$ and $\sin(n \omega t)$, each of which in turn is a finite polynomial combination of $\cos(\omega t)$ and $\sin(\omega t)$. The overall order of the dynamical system necessarily increases but autonomous polynomial dynamics are still sufficient to capture the systems' dynamics; see \cite{Mcmillan20} for further details.
\par There are a few requisite remarks that need to be made in view of this auxillary function method. The first is that this auxillary function method computes the maximal long time average of trajectories across \textit{all} initial conditions, but it does not give the knowledge of \textit{which} initial conditions or trajectories attain the maximal long time average. If one were interested in the analysis of the particular solution that achieves the maximal long-time average, other methods would need to be developed. This is perhaps one potential drawback of the method. However, in many applications, the knowledge of these extremal averages is sufficient for the study at hand. This method is particularly well-suited for systems subject to random initial conditions and when long-time statistics are of the utmost importance. Finally, we note for the interested reader that there are extensions of this auxillary function method to partial differential equations applications; see the work of \cite{Goluskin19, Rosa20, Goulart12} for explication.
\section{Model Setup}
Our general model of interest is a parametrically driven, coupled oscillator system of the form
\begin{equation}
    \begin{split}
        &\ddot{x}_1+{\omega_0}^2[1+h\sin(\gamma t+\frac{\phi}{2})]x_1+\omega_0g \dot{x}_1-\omega_0r\dot{x}_2=0\\
        &\ddot{x}_2+{\omega_0}^2[1+h\sin(\gamma t-\frac{\phi}{2})]x_2+\omega_0g\dot{x}_2+\omega_0r\dot{x}_1=0.
        \label{linear_osc}
    \end{split}
\end{equation}
In Equation (\ref{linear_osc}), $\omega_0$ denotes the proper frequency of the oscillators and $g$ is the intrinsic loss term, which is taken to be equal for both oscillators for the sake of simplicity. The $h$ and $\gamma$ terms are the intensity and frequency of the parametric stiffness terms, respectively. We focus on $h \in [0,1]$ and $\gamma \in [-3,3]$, as this is the regime which encapsulates most physically relevant phenomena. The coupling strength $r$ describes an energy preserving coupling between the oscillators, which corresponds to rotations in the $(x_1,x_2)$ plane and preserves the system's total energy. We note that in the limit as $r\rightarrow 0$ that the system in Equation (\ref{linear_osc}) becomes equivalent to two, decoupled parametric oscillators described by two, damped Mathieu equations. Finally, $\phi$ denotes the phase difference between the oscillators. Just as in \cite{Marcello19}, we will focus on the cases $\phi=0,\frac{\pi}{2},\text{and}\, \pi$, and we remark that the system in Equation (\ref{linear_osc}) for the aforementioned $\phi$ values exhibits three potential resonance frequencies at $\gamma=2\Omega_r$ and $\gamma=2\Omega_r\pm \omega_0 r$, where $\Omega_r=\omega_0\sqrt{1+(r^2-g^2)/4}$.
\par For the purposes of the auxillary function method for long-time averages, the system in Equation (\ref{linear_osc}) can be written as a first order, coupled system with exclusively polynomial terms:
\begin{equation}
    \begin{split}
        \dot{x}_1&=y_1\\
        \dot{x}_2&=y_2\\
        \dot{y}_1&=-\omega_0^2[1+hx_3]x_1-\omega_0g y_1+\omega_0 r y_2\\
        \dot{y}_2&=-\omega_0^2[1+hx_3]x_2-\omega_0g y_2-\omega_0 r y_1\\
        \dot{x_3}&=\gamma x_4\\
        \dot{x}_4&=-\gamma x_3\\
        &\text{Subject to: $x_3^2+x_4^2=1$},
        \label{poly}
    \end{split}
\end{equation}
where $\phi=0$. The system in Equation (\ref{poly}) can be written similarly for additional values of $\phi$ via elementary trigonometric identities. The constraint $x_3^2+x_4^2=1$ is enforced via the S-procedure and forces the variables $x_3$ and $x_4$ to be  uniquely determined. However, we remark that the S-procedure is not enforced when $\gamma=\phi=0$, as Equation (\ref{poly}) becomes independent of trigonometric terms. We employ the auxillary function method for long-time averages on the quantity of interest 
\begin{equation}
    \Phi=x_1^2+x_2^2.
    \label{phi}
\end{equation}
Hence, computing $\overline{\Phi}^*$ gives the maximal long-time average of the summed, mean squared amplitudes, which allows us to distinguish between small and large amplitude, parameter dependent solutions. The intuition being that small or large amplitude solutions should generically correlate with asymptotically stable or unstable solutions, respectively. Since the maximal long-time averages will be either zero or infinite in the regions of stability and instability, respectively, the auxillary function method will allow us to trace out the threshold between stability and instability. The fineness of the  threshold, moreover, only depends on one's mesh size for the system parameters \footnote{The mesh size is 80 by 80 for every figure in this paper.}. With $\Phi$ as defined in Equation (\ref{poly}), the sum of squares program to solve becomes: 
\begin{equation}
\begin{gathered}
    \min \,\,U\\
    \text{s.t.}\quad U-\Phi -\textbf{f}(x_1,x_2,y_1,y_2,x_3,x_4)\cdot \nabla V+\\
    S(x_1,x_2,y_1,y_2,x_3,x_4)(1-x_3^2-x_4^2) \in \mathcal{S}[x_1,x_2,y_1,y_2,x_3,x_4]\\
     S(x_1,x_2,y_1,y_2,x_3,x_4)(1-x_3^2-x_4^2)  \in \mathcal{S}[x_1,x_2,y_1,y_2,x_3,x_4], 
    \end{gathered}
\end{equation}
where $\textbf{f}(x_1,x_2,y_1,y_2,x_3,x_4)$ is the polynomial vector field as defined in Equation (\ref{linear_osc}), $V$ is an auxillary function polynomial whose degree may vary, and we have enforced the S-procedure at all $\gamma$ values other than zero. The convex optimization problems in this paper were solved using Mosek\cite{Mosek} paired with Yalmip\cite{Yalmip}. 
\par We remark that similar computations using the auxillary function method have already appeared in \cite{Mcmillan20}, where the model of interest was the Duffing oscillator. In particular, it was shown that the auxillary function method accurately reproduces the Duffing equation's frequency response curve and parameter dependent hysteresis phenomena as derived by harmonic balance down to numerical precision. This further validates the auxillary function method as a tool to study oscillator dynamics.
\par To the author's knowledge, this is the first time that the auxillary function method has been applied to studying nonlinear, coupled oscillators or has been used to explicitly study dynamical stability. Moreover, as we suspect and as we will show, the dynamics of Equation (\ref{linear_osc}) are substantially more complicated than that which can be described by asymptotic methods. The reason is that in contrast to the Duffing equation, our model has two \textit{coupled}, nonlinear oscillators and the non-autonomous forcing appears parametrically instead of externally. Hence, the synchronization of the two oscillators plays a crucial role in their stability and that synchronization pivotally depends on one's choice of system parameters.

\par In Section V. and in a similar fashion to \cite{Mcmillan20}, the parameter dependent computations of $\overline{\Phi}^*$ will be used to recover stability region boundaries and the validity of these results will be corroborated via asymptotic analysis. 
\section{Asymptotic Analysis}
We will corroborate the findings of the auxillary function method for long-time averages via comparison with the approximate solutions established by Floquet theory. This is a standard method which is discussed by many sources on asymptotic analysis; the specific reference used here is \cite{Miller}.
\par Consider $x(t)\in \mathbb{R}$ satisfying
\begin{equation}
    \ddot{x}(t) + F(t)x(t) = 0,
    \label{flq1}
\end{equation}
where $F(t)$ is a periodic function with period $T$. Equation \eqref{flq1} is a second order differential equation and thus has two linearly independent solutions $x_1(t)$ and $x_2(t)$.
\par Since $x_1(t)$ and $x_2(t)$ are linearly independent, they span the solution space of Equation \eqref{flq1} In particular, there must exist $\textbf{M}\in \mathbb{R}^{2\times 2}$ such that
\begin{equation}
    \begin{bmatrix}
    x_1(t+T) \\ x_2(t+T) \end{bmatrix} = \textbf{M} \begin{bmatrix} x_1(t) \\ x_2(t) 
    \end{bmatrix}, 
    \end{equation}
where $\textbf{M}$ is called the \emph{monodromy matrix} of $x_1(t)$ and $x_2(t)$. A monodromy matrix can be constructed for any pair of linearly independent solutions of Equation \eqref{flq1}. One can show that the determinant of the monodromy matrix is $1$ and eigenvalues of the monodromy matrix are independent of the pair of solutions that one chooses to start with.

Suppose $\lambda_1$, $\lambda_2$ are eigenvalues of a monodromy matrix $\textbf{M}$ corresponding to solutions $x_1(t)$, $x_2(t)$. For each $\lambda_j$ there exists a linear combination $x_{\lambda_j}(t)$ of $x_1(t)$ and $x_2(t)$ such that
\begin{equation}
    x_{\lambda_j}(t+T) = \lambda_j x_{\lambda_j}(t).
\end{equation}
\emph{Floquet's theorem} states that we may write
\begin{equation}
    x_{\lambda_j} = e^{-i\mu_j t} u_{\lambda_j}(t),
    \label{flq2}
\end{equation}
where $u_{\lambda_j}(t)$ is periodic with period $T$ and $\mu_j \in \mathbb{C}$ is such that $\lambda_j = e^{-i\mu_j T}$. 
Since $u_{\lambda_j}(t)$ is periodic, one may expand $u_{\lambda_j}(t)$ as a Fourier series and recover the analytic solution $x_{\lambda_j}(t+T)$ to within arbitrary accuracy by recursively solving for the Fourier amplitudes.

\par  For Equation (\ref{linear_osc}) with $\phi = 0$, we generally follow the method outlined in \cite{Marcello19}, where the study focused on frequencies at or near parametric resonance due to an interest in phase-locking, but our novel contribution is that we do not limit our attention to resonant frequencies.
\par The equations in Equation (\ref{linear_osc}) can be decoupled by performing a change of basis and defining
\begin{equation}
    x_{\pm}(t) = x_1 \pm i x_2,
\end{equation}
Adding the first line of Equation \eqref{linear_osc} to $\pm i$ times the second line yields 
\begin{equation}
    \ddot{x}_{\pm} + \omega_0^2 \left(1 + h \sin(\gamma t) \right) x_{\pm} + \omega_0 \left(g\pm ir\right) \dot{x}_{\pm} = 0,
    \label{decoupled}
\end{equation}
where the solutions in the new basis $x_\pm(t)$ exhibit real and imaginary loss terms.
We then let

\begin{equation}
    x_{\pm}(t) = e^{-\frac{(g \pm ir) \omega_0 t}{2}} y_\pm(t)
    \label{decoupled_sol}
\end{equation}
for $y_\pm(t)$ to be determined.  Upon substituting Equation (\ref{decoupled_sol}) into the decoupled system of Equation (\ref{decoupled}), we arrive at:
\begin{equation}
    \ddot{y}_{\pm} + \omega_0^2 \left[1 - \frac{(g \pm ir)^2}{4}+h \sin(\gamma t) \right]y_{\pm}(t) = 0.
    \label{decoupledlinear}
\end{equation}
Since the coefficient of $y_{\pm}$ in Equation \eqref{decoupledlinear} is periodic with period $T = \frac{2 \pi}{\gamma}$, then according to Floquet Theory, we may write
\begin{equation}
    y_{\pm}(t) = e^{-i\mu t}f_\pm(t),
\end{equation}
where $f_\pm(t)$ are some periodic functions with period $T = \frac{2 \pi}{\gamma}$. Wherever $\text{Im}(\mu) < 0$, $y_+(t)$ will be unstable as $t\rightarrow \infty$. Since $f_\pm(t)$ are periodic, we may express them as Fourier series
\begin{equation}
   f_\pm(t) = \sum_{n \in \mathbb{N}} A_n^\pm e^{i n \gamma t}, 
   \label{fourier}
\end{equation}
where $A_n^\pm$ denotes the amplitude of the nth Fourier coefficient for $f_\pm(t)$, respectively. 
If we substitute Equation \eqref{fourier} into Equation \eqref{decoupledlinear} and collect the coefficients of $e^{in\gamma t}$, we get the recursion relation for the Fourier coefficients $A_n^\pm$ as derived by \cite{Marcello19}:
\begin{equation}
\begin{gathered}
    D_{\pm,n}(\mu) A_n^\pm +i\frac{\omega_0 h^2}{2}(A_{n+1}^\pm-A_{n-1}^\pm)=0,
    \label{recursion}
    \end{gathered}
\end{equation}
where $D_n(\mu)=\omega_0^2-(n\gamma-\mu^2)+i\omega_0(g+ir)(n\gamma-\mu)$. We note that if the driving intensity $h$ is less than one, $A_n^\pm$ is coupled most strongly to $A_{n\pm 1}^\pm$, with coupling proportional to $h^2$. Without loss of generality, we will consider $n=0$.
Therefore, one recovers a matrix equation of the form:
\begin{equation}
    \text{Det}\begin{pmatrix}
    \textbf{D}_{-1}(\mu) & \textbf{M}(h) & 0 \\
    \textbf{M}^*(h) & \textbf{D}_{0}(\mu) & \textbf{M}(h)\\
    0 & -\textbf{M}^*(h) & \textbf{D}_{+1}(\mu)
    \end{pmatrix}=0,
    \label{matrix2}
\end{equation}
where one defines
\begin{equation}
\begin{split}
    \textbf{D}_n(\mu)&=\begin{pmatrix}
    D_{+,n}(\mu) & 0 \\
    0 & \textbf{D}_{-,n}(\mu) 
    \end{pmatrix}
    \\
    \textbf{M}(h)&=\begin{pmatrix}
    \frac{ih\omega_0^2}{2} & 0 \\
    0 & \frac{ih\omega_0^2}{2} 
    \end{pmatrix}.
    \end{split}
\end{equation}
Note that in Equation (\ref{matrix2}) that enforcing the determinant to vanish ensures that there are non-trivial solutions to Equation (\ref{decoupledlinear}). In contrast, it was also shown in \cite{Marcello19} that the simplifying assumptions of $r=g=0$ and ignoring the coupling effects of $D_{+1}$ yield the matrix equation:
\begin{equation}
    \text{Det}\begin{pmatrix}
    \omega_0^2-(\gamma-\mu)^2 & \frac{i h \omega_0^2}{2} \\
    -\frac{i h \omega_0^2}{2} & \omega_0^2-\mu^2 \\
    \end{pmatrix}=0,
    \label{matrix3}
\end{equation}
The matrix equation in Equation (\ref{matrix3}) is a consequence of assuming the coupling effects as well as the higher order effects are negligible in the asymptotic analysis. Hence, we hereafter refer to the solutions of Equation (\ref{matrix3}), which is the same as presented in \cite{Marcello19}, as the simplified, uncoupled solution.
\par Similar to \cite{Marcello19}, we focus on the solutions of Equation (\ref{matrix2}) and Equation (\ref{matrix3}) that display parametric resonance at $\gamma=2\Omega_r, 2\Omega_r\pm \omega_0 r$; we present the solutions in the Results section. We show that the simplifying assumptions have highly non-trivial and significant effects on the system's regions of stability at and away from parametric resonance frequencies.
\section{Results}
We first compare the stability boundary of the results established via Equation (\ref{matrix2}) with the stability boundary of the simplified, uncoupled solution of Equation (\ref{matrix3}). We then compare the stability boundaries as predicted by asymptotic analysis with the stability boundary as predicted by the auxillary function method for long time averages.  
\subsection{Influence of Higher Order Effects}
\par Indeed, we employ MATLAB's Solve function to find the roots in $\mu$ of Equation (\ref{matrix2}), whose solutions are sufficiently complicated and long to not be included in the text \footnote{The expressions for the roots can be found on the lead author's web-page: sites.google.com/andrewmcmillan/Research.}; This yields six solutions. In Figure \ref{Fig1} \footnote{The authors remark that implicitly solving for the roots of $\mu$ is computationally sensitive to very small numerical error. Hence, the stability boundary seen in Figure \ref{Fig1} had miniscule variations depending on the software used. However, the general shape of the stability boundary is consistent across software.}, we plot the contours of the imaginary part of one of the combined boundaries of the six solutions of our more generalized result in comparison to the simplified, uncoupled solution to Equation (\ref{matrix3}) for a wide range of parameter values:
\begin{figure}[ht]
    \centering
    \includegraphics[width=7.5cm, height=6cm]{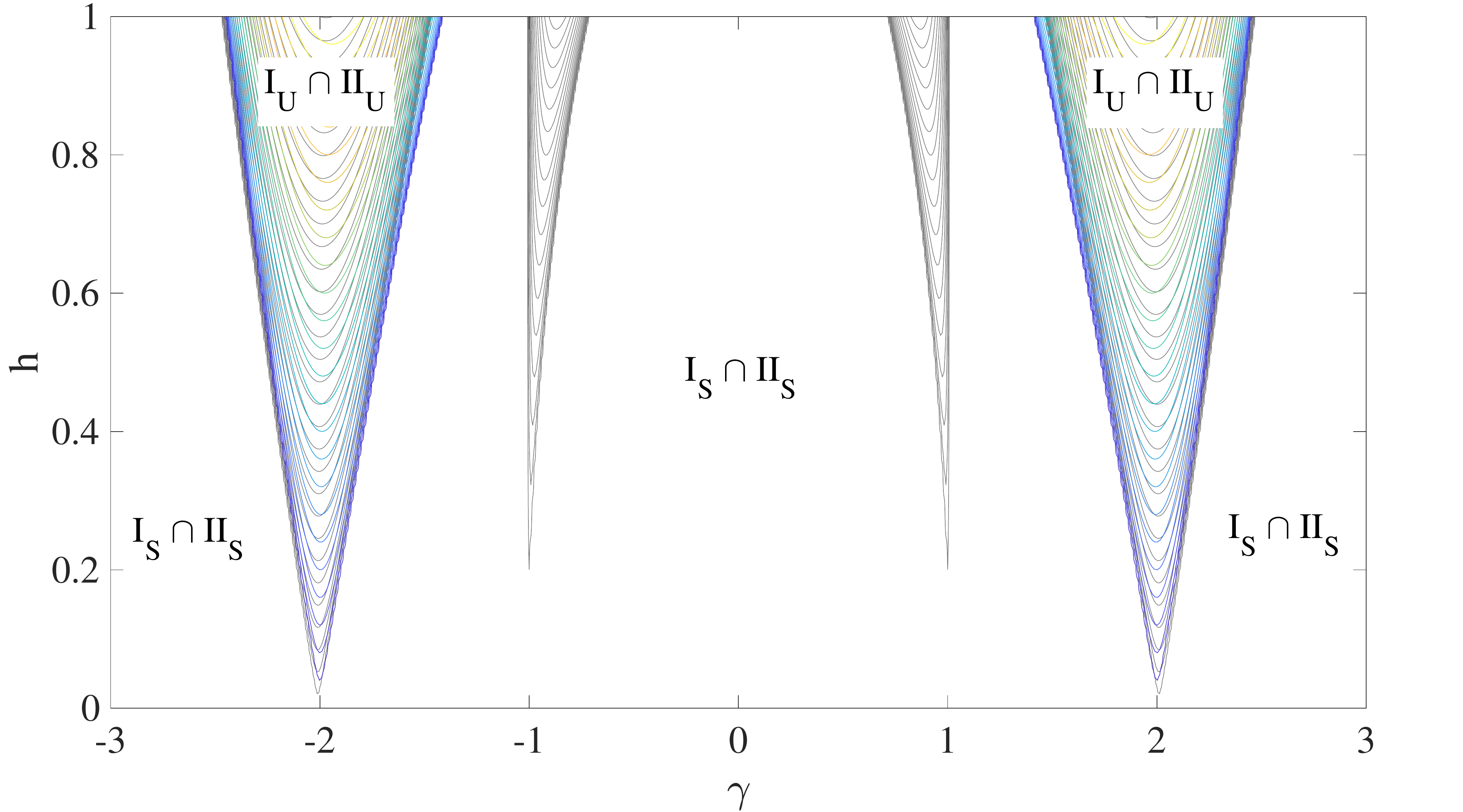}
    \caption{The simplified, uncoupled asymptotic solution to Equation (\ref{matrix3}) (shown with colored, solid contour lines) and the more general asymptotic solution to Equation (\ref{matrix2}) (shown with dashed, gray contour lines). The markers $\text{I}_\text{S}$, $\text{I}_\text{U}$, $\text{II}_\text{S}$, and $\text{II}_\text{U}$ denote the stable (S) or unstable (U) regions of the simplified (I) and general (II) solutions, respectively. The parameter values are $r=.2$, $g=.01$, and $\phi=0$.}
    \label{Fig1}
\end{figure}
\par In Figure \ref{Fig1}, we see that the more general, higher order asymptotic analysis solution is a dramatic and surprising improvement on the simplified, uncoupled solution, as the region of instability for the simplified, uncoupled asymptotic result is a subset of the region of instability for the more general asymptotic result. While there is good agreement between the two solutions for $|\gamma|>1$, there is a large region within the parameter space, $|\gamma|\leq 1$, for which the simplified, uncoupled asymptotic analysis solution predicts stability, while the more general, higher order solution reveals instability.
\par The fact that the simplified, uncoupled asymptotic expansion stability solution  fails so drastically for $|\gamma| \leq 1$ is quite a surprising finding. This means the simplified solutions are dangerously un-conservative in predicting regions of stability.
However, it then becomes natural to ask if the more general, higher order solution also fails to fully describe the stability boundary of the solutions to Equation (\ref{linear_osc}). That is, we may ask if including successively higher order terms in the asymptotic expansion would lead to such drastic changes in the stability boundary as seen between the two solutions in Figure \ref{Fig1}. In order to address this question, we compare the stability regions as predicted by the general asymptotic analysis with the stability regions predicted by the auxillary function method for long time averages.
\subsection{Asymptotic Analysis vs. The auxillary function method}
\par Upon computing the stability boundary with the auxillary function method for long time averages using the expression in Equation (\ref{phi}), we arrive at the solid black line seen Figure \ref{Fig2}:
\begin{figure}[h!]
    \centering
    \includegraphics[width=7.5cm, height=6cm]{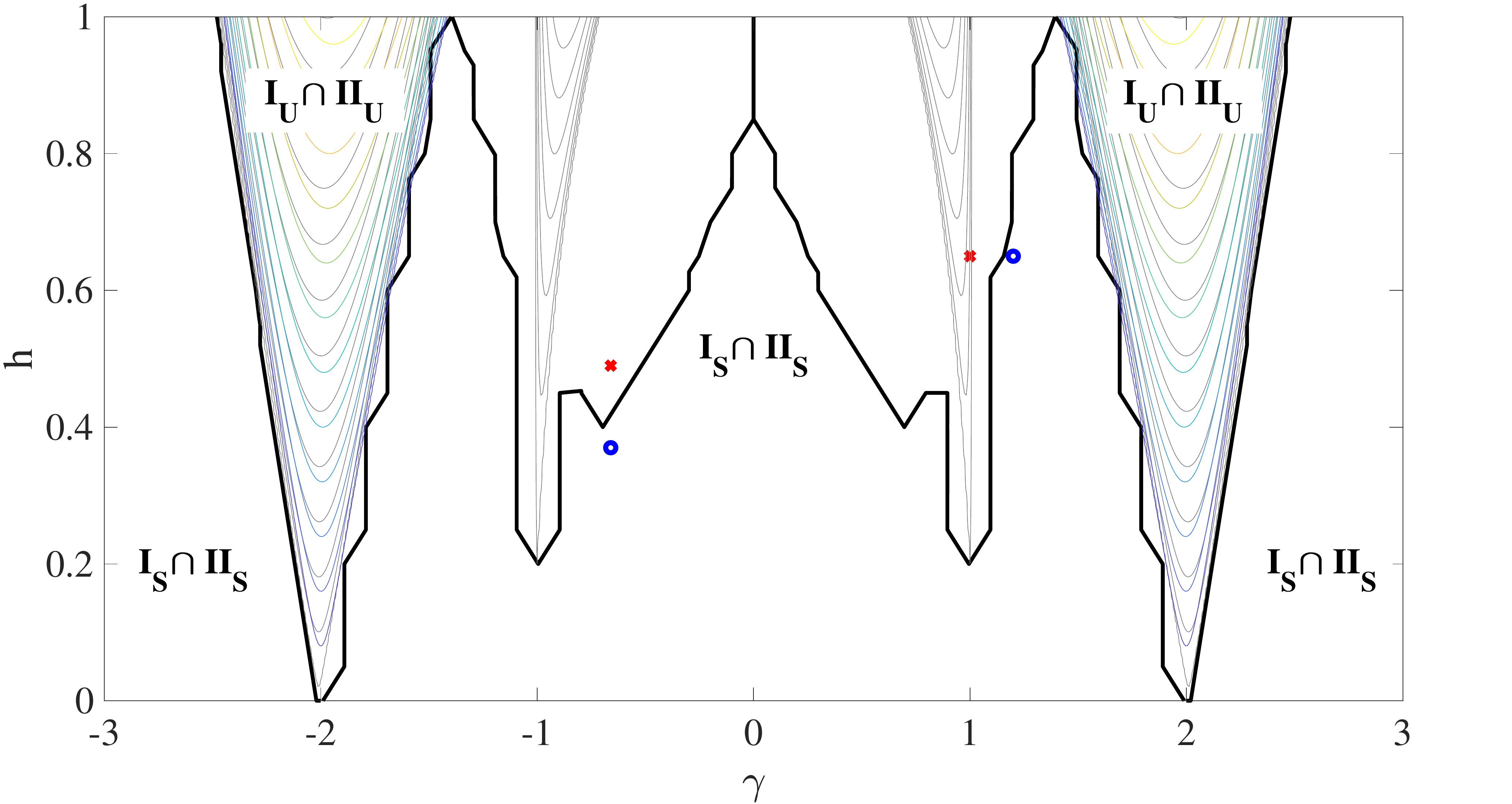}
    \caption{There is strong agreement between our auxillary function method (shown with a thick, solid black line) with the theoretical predictions of the simplified, uncoupled asymptotic analysis (shown with colored, solid contour lines) and the more general, higher order asymptotic analysis (shown with dashed, gray contour lines) for  $r=.2$, $g=.01$, and $\phi=0$. The open blue circles and red crosses indicate where direct numerical simulation results are shown in Figure \ref{Fig3} for stable and unstable points, respectively.}
    \label{Fig2}
\end{figure}
\par In Figure \ref{Fig2}, we see that the more general asymptotic analysis result agrees quite well with the results of the auxillary function method. In particular and most noticeably, both methods capture a pair of narrow, protruding tongues that occur at $\gamma\in \{-1,1\}$. However, there is still a quite large range of $\gamma$ values for which the auxillary function method predicts potential instability, but the asymptotic analysis solution to Equation (\ref{matrix2}) predicts stability. 
\par In order to validate the extended region of instability as predicted by the auxillary function method for long time averages results, we choose four points within the parameter space to perform  direct numerical simulation using ode45 in MATLAB.The results show that the auxillary function method for long time averages is corroborated via direct numerical simulation.  
\begin{figure}[ht]
    \centering
    \includegraphics[width=16cm, height=8.9cm]{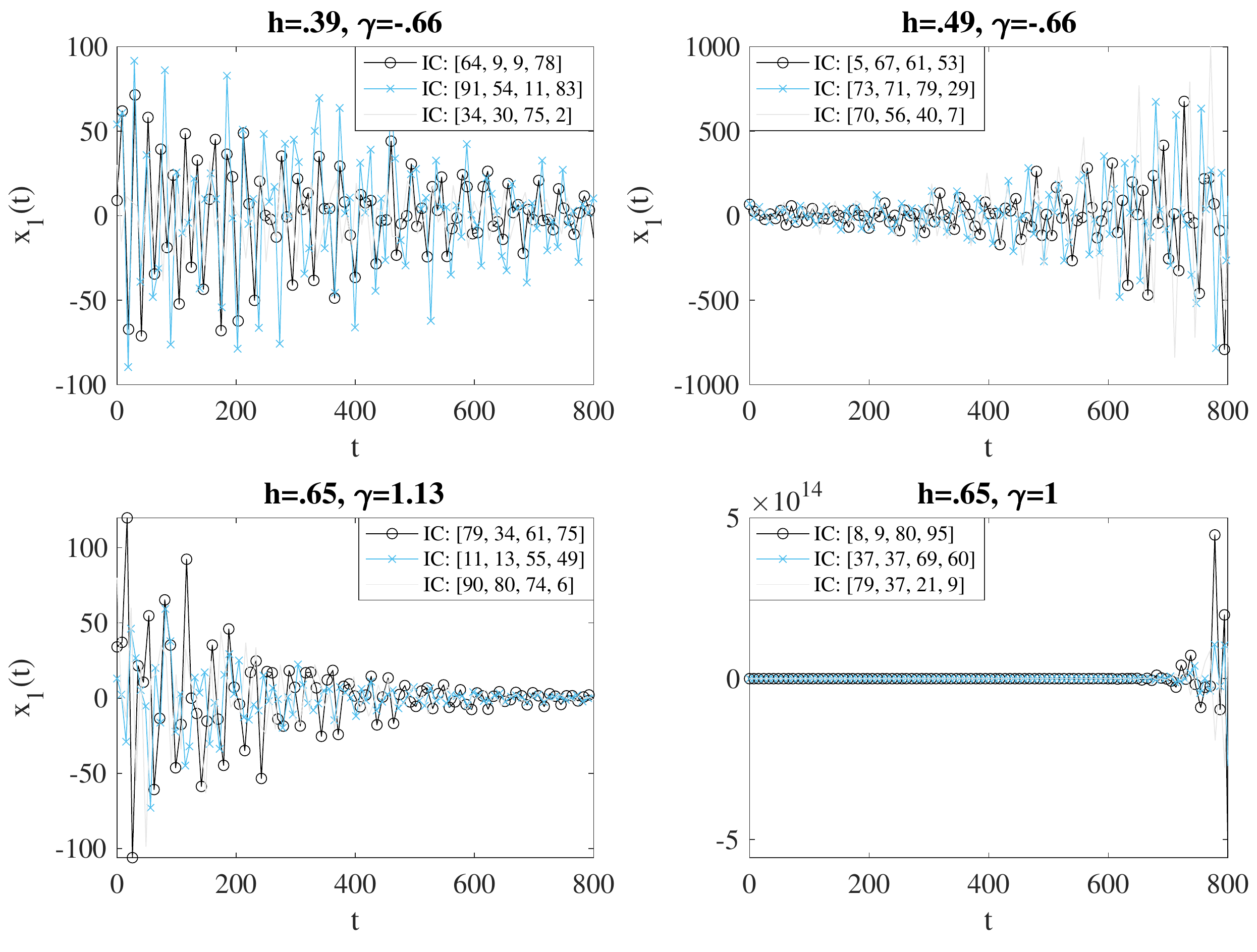}
    \caption{ The time histories predicted by DNS for the four points marked in Figure \ref{Fig2}. The left two plots correspond to the two blue, open circles just outside the region of instability for Figure \ref{Fig2}, and the two right plots correspond to the two red crosses just inside the region of instability for Figure \ref{Fig2}. In the above plot, an initial condition vector IC: [a,b,c,d] corresponds to $x_1(0)=a$, $x_2(0)=b$, $y_1(0)=c$, and $y_2(0)=d$ for Equation (17) with $r=0.2$, $g=0.01$, and $\phi=0$.}
    \label{Fig3}
\end{figure}

\par That is, for points marked by red crosses inside the region of instability, there is a blow up of the solutions, with the time histories shown on the right hand side of Figure \ref{Fig3}. However, for points marked by open, blue circles just outside the instability region experience decay, with the time histories shown on the left hand side of Figure \ref{Fig3}.
\par Additionally, Figure \ref{Fig3} displays three trajectories for three randomly chosen initial conditions for each of the four chosen points in parameter space. We also note that the above figure only plots the trajectory of $x_1(t)$ for our model, but when we plotted the trajectories for $x_2(t)$ a similar pattern of blow-up or stability still held except that $x_2(t)$ blew up at a rate much slower in comparison to $x_1(t)$. This is to be expected as the coupling term appears negatively in the equation for $x_1(t)$ and positively in the equation for $x_2(t)$. Hence, the $x_2(t)$ oscillator feeds energy into the $x_1(t)$ oscillator, so we would expect a faster blow-up for $x_1(t)$.
\par This is perhaps a case study example of how the auxillary function method can lend itself to finding regions of instability where perturbative methods fail.
Strikingly, the perturbative method for this system seems to fail quite drastically with a large portion of the stability diagram not being captured by the perturbative method but instead by the auxillary function method. Therefore, the auxillary function method has the advantage of being able to recover the true stability boundary both at and away from parametric resonance, while doing so in a computationally efficient way. We  remark that the auxillary function method is quite computationally efficient despite the polynomial representation, seen in Equation (\ref{poly}) , containing five degrees of freedom. 
In particular, the CPU times corresponding to the computation of $\overline{\phi^*}$ at the points $(\gamma, h)=(-1.13, .650), (-1.00, .650), (-.660, .390),$ and $(-.660, .490)$ are 3.4844, 2.6744, 2.8574, and 2.5737 seconds for a degree 8 auxillary function and a degree 6 S-procedure enforcement on nothing more than a standard laptop using a single core, 2.2 GHZ processor\footnote{The source code for the various points is also available on the lead author's web-page.}.
\par In light of the above results, it then becomes natural to analyze how varying the system parameters changes the above regions of stability. In particular, we study the effects of varying the coupling strength in the next section.

\subsection{Effects of Damping and Coupling Parameters}
\par In order to study the affects of varying the coupling strength, we consider Equation (21) for $r=.4$ and compare the stability boundary with the prior results of the both auxillary function method as well as the asymptotic analysis for $r=.2$ and $g=.01$. We do not consider negative values for $r$ because the coupling term appears negatively in the first equation and positively in the second equation of Equation (21); hence, a sign change of $r$ just swaps the role of $x_1(t)$ and $x_2(t)$, respectively. Also, we focus on $\gamma \in [-1.5,1.5]$, as this is where the discrepancies between the auxillary function method and asymptotic analysis were established in the previous section.

\par
\begin{figure}[h!]
    \centering
    \includegraphics[width=7.5cm, height=7cm]{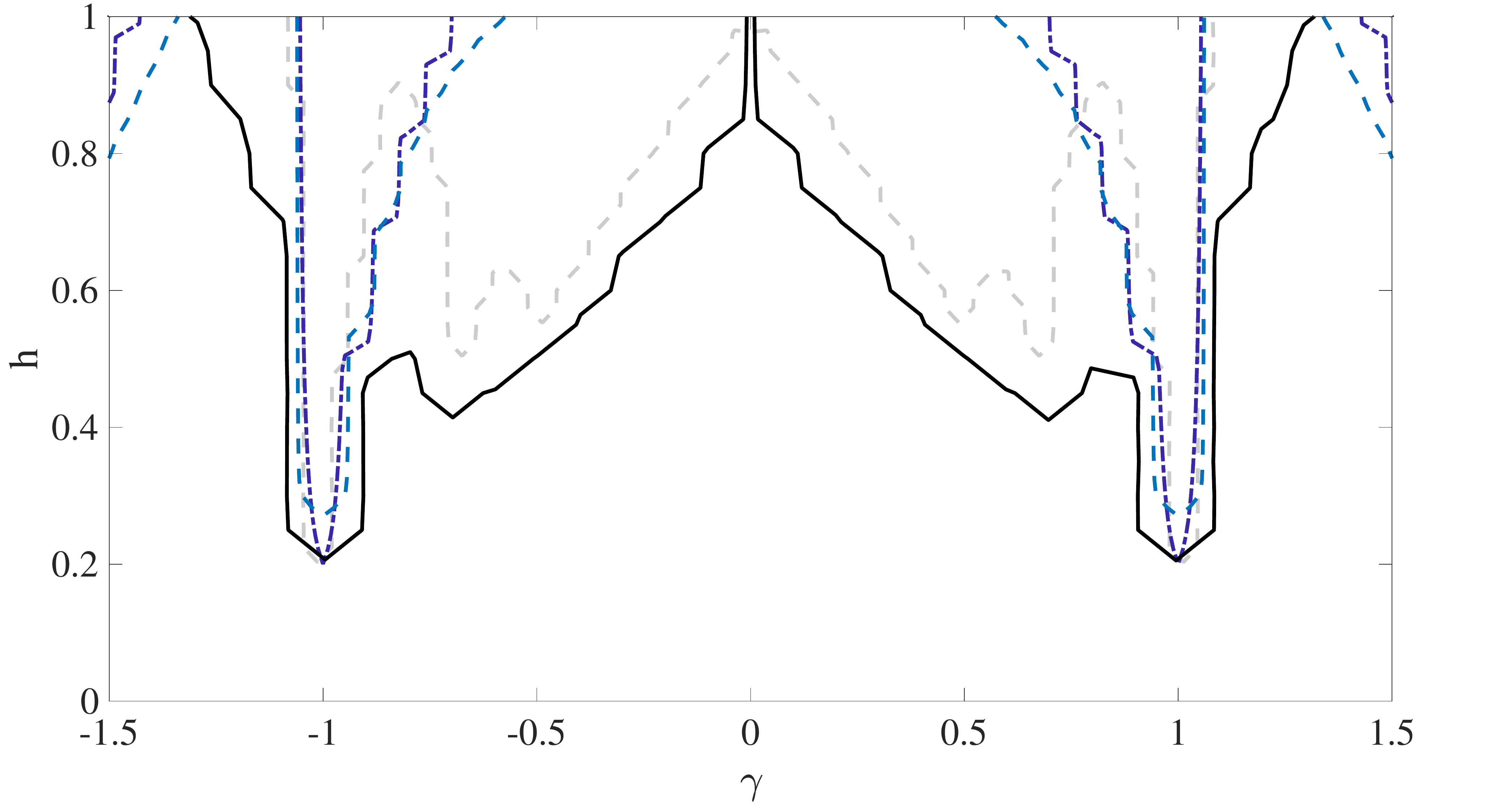}
    \caption{The stability boundary as predicted by the auxillary function method is shown with a thick, solid black line and a thick, dotted gray line for $r=.2$ and $r=.4$, respectively, while the asymptotic analysis boundary is shown with a thin, purple dashed line and a thin, blue dotted line for $r=.2$ and $r=.4$, respectively, with $g=.01$ and $\phi=0$. }
    \label{Fig4}
\end{figure}

 Performing the same procedure with the auxillary function method as in the previous section as well as incorporating the predictions of the asymptotic analysis, we find in Figure \ref{Fig4} that the region of stability for the equations with $r=.2$ is a subset for the region of the stability for the equations with $r=.4$. Intuitively, one would expect this to be true because as the r-value increases, the rate of energy transfer from one oscillator to another increases and hence initial transients are far less likely to become unstable.
\par Figure \ref{Fig4} shows that the discrepancies between the auxillary function method and the asymptotic analysis solution for $r=.4$ persist. Moreover, the discrepancies seem to grow as the value of $r$ is increased. Therefore, as the coupling strength between the oscillators increases, the asymptotic analysis solution to Equation (\ref{matrix2}) increasingly fails to capture the true stability boundary of the system, because the coupling terms amplify the relative importance of the higher order terms, particularly away from parametric resonance as the sine and cosine terms do not cancel. We also note that Figure \ref{Fig4} is for $g=.01$, but if one is interested in the affects of $g$, we can similarly compute regions of stability for fixed $r$ and different g values. The expected result would simply be a shift up of the region of stability appearing in Figure 1. The reason is that $g$ appears on the damping term and hence increasing the damping coefficient should increase the size of the region of stability. 
\subsection{Effects of the Parametric Phase}
Finally, we study the effects of varying the phase value $\phi$. In a fashion similar to \cite{Marcello19}, we focus on $\phi=0,\frac{\pi}{2},\pi$. However, it is important to remark that \cite{Marcello19} was unable to determine a closed form, asymptotic solution for varying $\phi$. We also could not find a closed form solution to the analog of Equation \ref{matrix2} with $\phi\neq 0$. The ability to find the stability boundary for general system parameters proves to be another benefit of the auxillary function method. Indeed, performing the analogous computations via the auxillary function method and choosing the same $\Phi$ that appears in Equation (13), we arrive at Figure \ref{Fig5}:
\begin{figure}[h!]
    \centering
    \includegraphics[width=7.5cm, height=7cm]{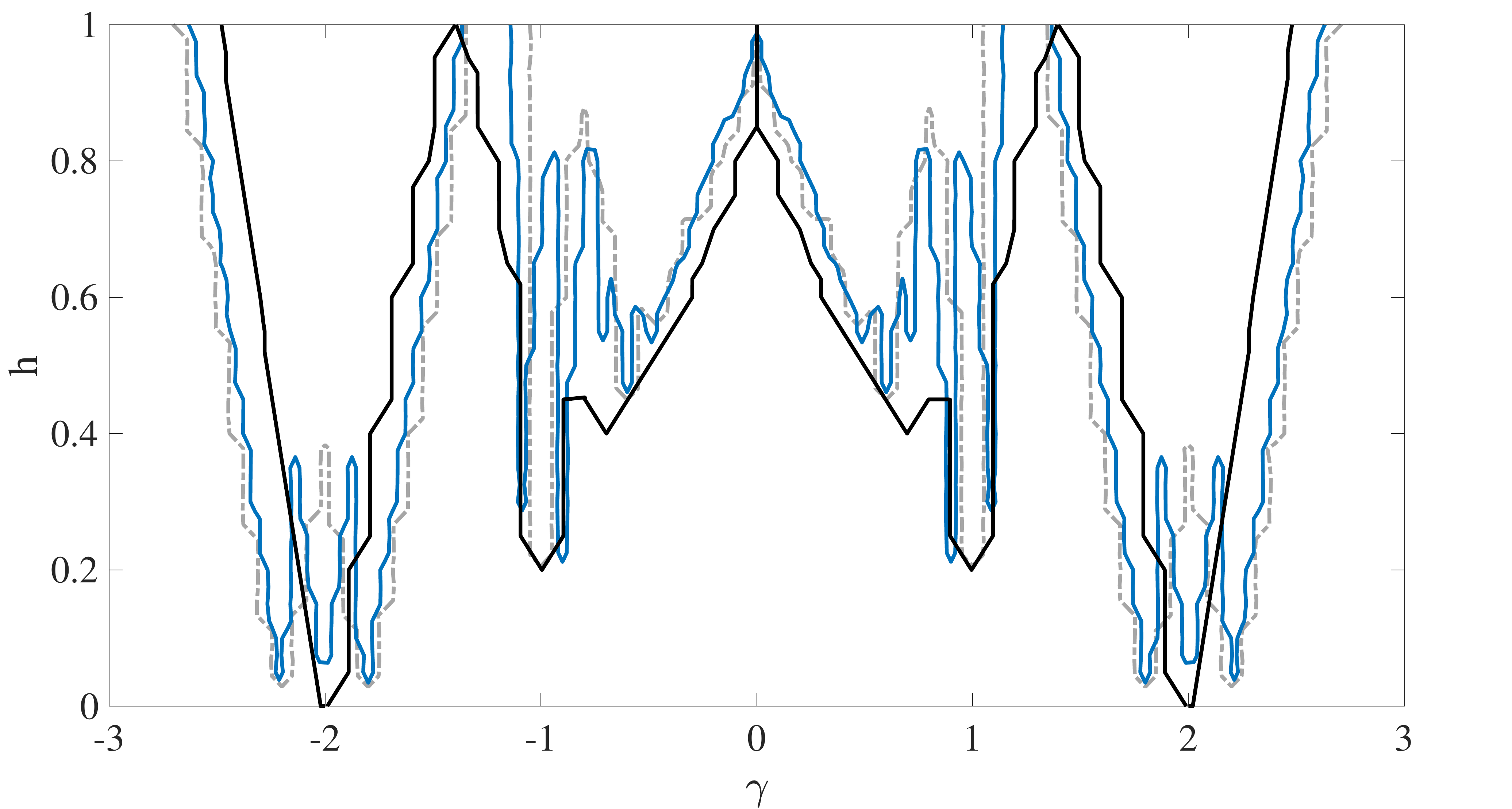}
    \caption{The stability boundary as predicted by the auxillary function method for $\phi=0$ (shown with a thick, solid black line), $\phi=\frac{\pi}{2}$ (shown with a thin, solid blue line), and  $\phi=\pi$ (shown with a thin, dashed gray line).}
    \label{Fig5}
\end{figure}
\\
Figure \ref{Fig5} shows that the auxillary function method reproduces the protruding tongues around $\gamma=\pm 2\Omega_r$ and $\gamma=\pm 2\Omega_r\pm \omega_0r$ with $\phi=0$, $\phi=\pi$, and $\phi=\frac{\pi}{2}$ corresponding to one, two, and three tongues respectively. Moreover, note that varying $\phi$ has non-trivial effects throughout the entirety of parameter space, and not only about $\gamma=\pm 2\Omega_r, \pm 2\Omega_r\pm \omega_0 r$, as varying the phase of the parametric oscillator term can lead to more localized tongue formation. 
\par Since the phase of a parametric oscillator can often be random, it is important to identify the instability boundary when varying phases are taken into consideration, especially when there is non-trivial disagreement between the true stability boundary and the boundary as predicted by asymptotic analysis.

\section{Conclusion and Discussion}
\par In summary, we have investigated the higher order effects caused by coupling parameters on the stability boundary of a non-linear, parametrically driven coupled oscillator system across a broad range of modulation frequencies. We have shown that the stability boundary of a non-linear, parametrically driven oscillator system as predicted by simplified, second order asymptotic methods ignoring the coupling terms between the oscillators can differ quite substantially at modulation frequencies away from parametric resonance frequencies. The simplified, asymptotic solution is un-conservative compared to a full, second order asymptotic solution when coupling terms are not ignored, which can in turn differ and still be non-trivially un-conservative in comparison to the true stability boundary at modulation frequencies away from the parametric resonance frequencies. The differences are caused by both neglecting the coupling terms and higher order effects. 
\par This is a primary drawback of critical importance for asymptotic methods. That is, the validity of the asymptotic results depends pivotally on one's choice of the approximations. However, there is currently no way to know which order will be sufficient to capture the true instability boundary, which depends on the strength of the coupling terms, $g$ and $r$, as well as the phase of the parametric oscillator term $\phi$. The solution is expected to be even more sensitive to these system parameters if the resonance frequency of the individual oscillations are different or if additional non-linearities are present. 
\par Hence,  we have shown that the auxillary function method for long-time averages is an efficient and robust means of computing the true long-time averages and true regions of stability across all possible initial conditions without the need of ad-hoc approximations. Moreover, this auxillary function method has the advantage of being able to compute regions of stability both at and away from parametric resonance.
\par The differences between the true stability boundary and the approximate stability boundary are immaterial if one is operating withing a region of parameter space where both boundaries agree. However, if one is operating within a region of parameter space for which they disagree, this may be quite problematic for both experimental and real-world implementation- especially without the knowledge of operating within one of these regions of disagreement. This point is exacerbated by the fact that, at least in studying Equation (\ref{linear_osc}), we discovered two, very narrow protruding tongues in parameter space for which the system was potentially unstable. Moreover, these tongues occur at very naturally occuring, non-pathological values of the driving frequency $\gamma$.
\par In the context of applications, our results have profound implications on the reliance of asymptotic methods across a variety of applied sciences. If one was interested in establishing a dynamical system via data-driven methods, the instability of these very narrow and protruding regions of parameter space could be completely missed, especially if data sampling was on the order of the size of the narrow tongues. If one was interested in studying synchrony behavior or global phase locking, asymptotic methods may fail to capture sensitivity to a system's parameters.
\par Hence, it is recommended that the instability region be determined and considered when planning new experiments or full-scale implementation, as disregarding the stability boundary or only relying on ad-hoc approximations may have costly consequences.
\section*{Dedication}
We are greatly indebted to Dr. Charles Doering (1956- 2021)- a great mentor, friend, and colleague. His insights, commentary and suggestions immensely impacted this work, and the entire scientific community will continue to miss his absence. 
\section*{Acknowledgments}
Support for this research is provided by the U.S. Office of Naval Research, contract No. N00014-18-C-1025 managed by Ms. Deborah Nalchajian.
\section*{Author Declarations} 
The authors acknowledge no personal, professional, or financial conflicts of interest.

\bibliographystyle{unsrt}

\end{document}